\documentclass[12]{article}

\usepackage{amssymb,latexsym,amsmath}     % Standard packages
\usepackage{dsfont}
\usepackage{upgreek}
\usepackage{graphics}
\usepackage{graphicx}
\usepackage{lscape}
\usepackage{multirow}
\usepackage{amssymb}
\usepackage{dsfont}
\usepackage{authblk}
\usepackage[switch,columnwise]{lineno}
\begin{document}
\title{A new 3-parameter extension of generalized lindley distribution}

\title{A new 3-parameter extension of generalized lindley distribution}
\author[1]{Deepesh Bhati}
\author[2]{Mohd. Aamir Malik}
\author[3]{K K Jose} 
\affil[1] {Department of Statistics, Central University of Rajasthan, India.  deepesh.bhati@curaj.ac.in}
\affil[2] {AU Financer(India) Ltd., Jaipur, Rajasthan, India, aamirmalik.stats@gmail.com}
\affil[3] {kkjstc@gmail.com}
\maketitle

\begin{abstract}
Here, we introduce a new class of Lindley generated distributions which results in more flexible model with increasing failure rate (IFR), decreasing failure rate(DFR) and upside down hazard functions for different choices of parametric values. We explore, various distributional properties including limiting distribution of extreme order statistics explored. Maximum likelihood estimators and the confidence intervals of the parameters are obtained. The applicability of the proposed distribution is shown through modelling two sets of real data on bladder cancer patients and waiting time in a queue. Further, we carry out stress-strength analysis for applying the model in system reliability studies.
\end{abstract}



\section{Introduction}
\indent Modelling of lifetimes of materials, organisms,devices plays an important role in biological and engineering services. Recently, a member of lifetime distributions have been introduced by various authors Nadarajah et al.(2012), Ristic(2012), Jose et al.(2014). They help in the modelling of length of life length data from various contents.

Lifetime distribution are used to describe, statistically, the length of the life of a system, a device, and in general, time-to-event data. Lifetime distributions are frequently used in fields like reliability, medicine, biology, engineering, insurance etc. Many parametric models such as exponential, gamma, Weibull have been frequently used in statistical literature to analyse lifetime data.

\indent Recently, one parameter Lindley distribution has attracted researchers for its potential in modelling lifetime data, and it has been observed that this distribution has performed excellently well in many applications. The Lindley distribution was originally proposed by Lindley\cite{r19} in the context of Bayesian statistics, as a counter example to fiducial statistics. The distribution can also be derived as a mixture of exp($\theta$) and {gamma(2, $\theta$)}. More details on the Lindley distribution can be found in Ghitany et al.(2008). Nastic et al.(2015) developed auto--generated time series models with Lindley marginal distribution and applied it to model time series data.

Ghitany et al(2013) introduced a power Lindley distribution and carried out associated inferences. Azgharzadeh et al.(2013) introduced a new lifetime distribution by compounding Poisson-Lindley distribution. Liyanage and Pararai(2014) introduced an exponential power Lindley distribution and explore its properties.

\indent A random variable $X$ is said to have the Lindley distribution with parameter $\theta$ if its probability density is defined as:
\begin{equation} \label{1}
r_X(x;\theta)=\frac{\theta^2}{(\theta+1)}(1+x) e^{-\theta x} ;x > 0,\theta>0.
\end{equation}
The corresponding cumulative distribution function is
\begin{equation*}
R(x)=1-\frac{e^{-\theta x}(1+\theta+\theta x)}{1+\theta}; x \in \mathbb{R}^+, \theta>0,
\end{equation*}
\indent Alzafreh et al (2013,2014) introduced a new class of distributions called T-X family. Ghitany et al.(2011) have introduced a two-parameter weighted Lindley distribution and have pointed out its usefulness, in particular, in modelling biological data from mortality studies.
Bakouch et al.(2012) have introduced extended Lindley (LE) distribution; moreover and Ebatal(2014) introduced a transmuted Lindley-geometric distribution and by transmuting and compounding Lindley and geometric distributions. Adamidis and Loukas(1998) have introduced a new lifetime distribution with decreasing failure rate. Shanker et al.(2013) have introduced a two-parameter Lindley distribution. Zakerzadeh et al.(1998) have proposed a new two parameter lifetime distribution and studied its properties. Hassan(2014) has introduced convolution of Lindley distribution. Ghitany et al.(2015) worked on the estimation of the reliability of a stress-strength system from power Lindley distribution. Elbatal et al.(2013) has proposed a new generalized Lindley distribution by considering the mixture of two gamma distributions.Ali(2013) considered the mean residual life function and carried out stress-strength analysis under different loss functions for Lindley distribution in the counted of Bayesian Inference.\\
Zografos and Balakrishnan (2009),proposed a kind of gamma-generated family. Based on a baseline continuous
distribution $F(x)$ with survival function $\bar{F}(x)$ and density $f(x)$, they defined the cumulative distribution function (cdf) and probability density function (pdf) as

\begin{equation*}
G_X(x)= \frac{1}{\Gamma(\alpha)}\int\limits_{0}^{-\log\left(\bar{F}(x) \right)} t^{\alpha-1}e^t dt, \quad \quad x \in \mathds{R}^+, \alpha > 0,
\end{equation*}\\
\\
The corresponding pdf is obtained as:
\begin{equation*}
g(x)= \frac{1}{\Gamma(\alpha)} \left(-\log \bar{F}(x)\right) \bar{F}(x)^{\alpha-1} f(x), \quad \quad x \in \mathds{R}^+, \alpha > 0,
\end{equation*}\\

\section{Distributional Properties of EGL Distribution}

We consider a new family of distribution generated by an integral transform of the pdf of a random variable $T$ which follows one parameter Lindley distribution. The distribution function of this new family is given as:
\begin{equation*}
G_X(x)= \int\limits_{0}^{-\log(\bar{F}(x))} r(t) dt, \quad \quad x\in \mathds{R}^+, \theta > 0
\end{equation*}
substituting $r(t)$ from equation (\ref{1}), we get the new family of distribution with d.f.
\begin{equation} \label{2}
G_X(x)= \frac{\theta^2}{1+\theta}\int\limits_{0}^{-\log(\bar{F}(x))} (1+t)e^{-\theta t} dt, \quad \quad x\in \mathds{R}^+, \theta > 0
\end{equation}
where $ \theta > 0 $ and ,the corresponding probability density function (pdf) is given by
\begin{equation} \label{3}
g(x)=\frac{\theta^2}{1+\theta}\left(1-\log \bar{F}(x)\right) \bar{F}(x)^{\theta-1} f(x),
\end{equation}
In this formation, we consider $F(x)$ corresponding to Extended exponential distribution with survival function $e^{1-(1+\lambda x)^\alpha}$ which yields the distribution function of the new distribution as 
\begin{equation} \label{4}
G(x)=1-\frac{e^{\theta -\theta  (1+x \lambda )^{\alpha }}\left(1+\theta  (1+x \lambda )^{\alpha }\right)}{1+\theta}
\end{equation}  
with corresponding density
\begin{equation}   \label{5}
g(x)=\frac{ \alpha  \theta ^2 \lambda  (1+x \lambda )^{2 \alpha -1}e^{\theta -\theta (1+x \lambda )^{\alpha }}}{1+\theta}.
\end{equation}
\\
We refer to random variable with survival function (\ref{4}) as Extended Generalized Lindley(EGL) distribution with parameters $\alpha$,$\theta$ and $\lambda$  and denote it by EGL($\lambda$,$\theta$,$\alpha $).\\ 

\noindent \textbf{Proposition 1:} If $X \sim \text{EGL}(\lambda,\alpha,\theta)$ then random variable $Y=(1+\lambda X)^{\alpha}-1$ follows $Lindley(\theta)$. \\

This paper is arranged as follows. In section 2 we introduce the new Lindley generated distribution and study its properties. Section 3 deals with maximum likelihood estimation. In section 4, we conduct the stress- and strength analysis estimation.. Section 5 is devoted for application on real data set.

\subsection{Shape of the density}
In this section, we introduce and study the distributional properties of the EGL. In particular, if $X \sim EGL(\lambda,\theta,\alpha)$ then the shapes of the density  and hazard function, moments, the density of the $r^{th}$ order statistics, and other important measures of the ELD are derived and studied in detail. \\

\noindent For the density function of the EGL distribution, the first and the second derivatives of $\log g(x)$ are
\begin{equation*}
\frac{d}{dx}\log g(x)=-\frac{\lambda  \left(1+\alpha  \left(-2+\theta  (1+x \lambda )^{\alpha }\right)\right)}{1+x \lambda }
\end{equation*}
and 
\begin{equation*}
\frac{d^2}{dx^2}\log g(x)=\frac{\lambda ^2 \left(1-2 \alpha -(-1+\alpha ) \alpha  \theta  (1+x \lambda )^{\alpha }\right)}{(1+x \lambda )^2}
\end{equation*}
Hence the mode of EGL distribution is given by the following theorem. \\

\noindent \textbf{Theorem 1:} The probability density function of EGL$(\lambda,\theta,\alpha)$ is unimodal for $\alpha>\frac{1}{2-\theta} \cap 0<\theta<2$ and is given as 
\begin{equation*}
x_0=-\frac{1}{\lambda}+\frac{1}{\lambda}\left(\frac{2 \alpha-1}{\alpha\theta}\right)^{\frac{1}{\alpha }}
\end{equation*}

\noindent \textit{Proof:} For $\alpha >1$, $\frac{d^2}{dx^2} \log g(x) <0$ i.e. the density function $g(x)$ is log-concave. Note $\left( \log g'\right)(\infty)=-\infty$  and $\left( \log g'\right)(\infty)=\lambda \left(\alpha(2-\theta)-1\right) >0 $ for $\alpha >\frac{1}{2-\theta}$ and $0<\theta<2$. This implies that for $\alpha >\frac{1}{2-\theta}$ and $0<\theta<2$, $g(x)$ has unique mode at $x_0$, where $x_0=-\frac{1}{\lambda}+\frac{1}{\lambda}\left(\frac{2 \alpha-1}{\alpha\theta}\right)^{\frac{1}{\alpha }}$, is the solution of the equation $\frac{d}{dx}\log g(x)=0$, for $\theta>2$, $\frac{d}{dx}\log g(x) <0$ i.e. $g(x)$ is decreasing in $x$.\\
\noindent Further $g(0)=\frac{\alpha \theta^2 \lambda}{1+\theta}=\alpha \lambda r(0)$, therefore $g(0) <(>) r(0)$ according as $\alpha \lambda <(>) 1$.
\\

%\begin{flushright}
%$ \Box $
%\end{flushright}
%\begin{figure} \label{f1}
%\centering
%\includegraphics[width=0.9\textwidth]{combined_PDF}
%\caption{PDF plot for various values of $\lambda$ and $\theta$.}
%\end{figure}

\noindent \textbf{Theorem 2:} The hazard function of EGL Distribution is decreasing, upside down and increasing according as $\left(0<\alpha<\frac{1}{2}\right)\bigcup\left(\frac{1}{2}<\alpha<1 \bigcap \theta > \frac{1-2\alpha}{\alpha(\alpha-1)}\right)$, $ \left(\frac{1}{2}<\alpha<1\bigcap \theta> \frac{1-2\alpha}{\alpha(\alpha-1)}\right) $ and $(\alpha>1) \cap (\theta>0)$ respectively.\\
\\
\noindent \textit{Proof:} Considering the hazard rate function (hrf) of the EGL distribution given by
\begin{equation}  \label{7}
h(x)=\frac{\alpha  \theta ^2 \lambda  (1+x \lambda )^{-1+2 \alpha }}{1+\theta  (1+x \lambda )^{\alpha }},
\end{equation}
\\
\noindent and using theorems of Glaser(1980), we can discuss the shape characteristics of the hrf of EGL Distribution. The function $ \eta (x) = -g^{'}(x)/g(x)$ for EGL Distribution is given by
\[
\eta (x)=\frac{\lambda  \left(1+\alpha  \left(-2+\theta  (1+x \lambda )^{\alpha }\right)\right)}{1+x \lambda } \quad \text{and} \quad \eta^{'}(x)=\frac{\lambda^2}{(1+x \lambda )^2}u(x)\
\]
\noindent where 
\[u(x)=\left(2\alpha-1 +(\alpha-1)\alpha \theta (1+\lambda x)^{\alpha}\right) \quad \text{and} \quad 
u'(x)=\lambda\alpha^2(\alpha-1)\theta(1+\lambda x)^{\alpha-1}\]

\noindent For $0<\alpha<\frac{1}{2}$, function $u(x)<0,\forall \lambda, \theta >0$ hence $\eta"(x)<0 \, \forall \, x$, hence from theorem(b) of Glaser(1980), hazard function is a decreasing function of $x$. Let us consider the case when $\frac{1}{2}<\alpha<1$, then $u(0)=2\alpha-1+\alpha(\alpha-1)\theta$ and $u(\infty)=-\infty$, if $\theta>\frac{1-2\alpha}{\alpha(\alpha-1)}$, then $u(x)<0$ implies $\eta'(x) <0 \, \forall \, x$, hence hazard function is decreasing, whereas, for $\frac{1}{2}<\alpha<1$ and $0<\theta \le \frac{1-2\alpha}{\alpha(\alpha-1)}$, $u'(x)<0 \, \forall \, x>0$ and $u(0)>0$, therefore $\exists$ a point $x_0$ such that $u(x)>0$ for $x \in (0,x_0)$ and $u(x)<0$ for $x \in \left[x_0,\infty\right)$ implies  $\eta'(x)>0$ for $x \in (0,x_0)$ and $\eta'(x)<0$ for $x \in \left[x_0,\infty\right)$. Hence from Glaser(1980) hazard function is upside down shape. Finally when $\alpha>1$, both $u(x), u'(x)$ are positive implies $u(x)$ as positively increasing function implies $\eta'(x)>0 \, \forall \,x>0$. Hence hazard function is increasing, which proves the theorem. \\
\noindent It can also be verified that
\begin{equation*}
h(0)=\alpha \lambda \left(\frac{\theta^2}{1+\theta}\right) \quad \text{and} \quad  \lim\limits_{x \rightarrow \infty} h(x) = \begin{cases}
0  \quad & \text{if} \quad 0<\alpha <1 \\
\alpha \lambda \theta \quad & \text{if} \quad \alpha =1 \\
\infty \quad & \text{if} \quad \alpha > 1 
\end{cases}
\end{equation*}
Hence from the above relation $h^{EGL}(0)=\alpha \lambda h^{L}(0)$. 

\noindent The pdf and hazard function for different parameter values are shown in the figure 1.

%\begin{figure} \label{f2}
%\centering
%\includegraphics[width=0.7\textwidth]{hazard}
%\caption{Hazard function plot for various values of $\lambda$ and $\theta$.}
%\end{figure}

\subsection{The Quantile Function of EGL distribution}
The cdf, $G_X(x)=1-\bar{G}(x)$, is given by using eq.(\ref{4}). Further, it can be noted that $G_X$ is continuous and strictly increasing so that the quantile function of $ X $ is $Q_X(\gamma)=G^{-1}_X(\gamma)$, $0<\gamma<1$. In the following theorem, we give an explicit expression for $Q_X$ in terms of the Lambert $ W $ function. For more details on Lambert $W$ function we refer the reader to Jo\'rda(2010) and also to  Nair et al.(2013) for discussion on quantile functions.
\\

\noindent \textbf{Theorem 3:} For any $\theta, \lambda, \alpha >0$, the quantile function of the EGL distribution is 
\begin{equation}  \label{8}
x_{\gamma}=G^{-1}(\gamma)=\frac{1}{\lambda}\left(-\frac{W_{-1}\left(e^{-\theta -1} (\theta +1) (\gamma-1)\right)+1}{\theta }\right)^{1/\alpha}-\frac{1}{\lambda },
\end{equation}
where $W_{-1}$ denotes the negative branch of the Lambert W function.\\ 
\\
\textit{Proof:} By assuming $p=(1+\lambda x)^\alpha$ the cdf can be written as 
\begin{equation*}
G_X(x)=1-\frac{e^{\theta(1-p)}(1+\theta p)}{1+\theta}
\end{equation*}
For fixed $\theta, \lambda, \alpha >0 $ and $\gamma \in (0,1)$, the $\gamma^{th}$ quantile function is obtained by solving $F_X(x)=\gamma$. By re-arranging the above, we obtain
\begin{equation*}
e^{\theta(1-p)}(1+\theta p)=(1-\gamma)(1+\theta)
\end{equation*}
It can be further written as
\begin{equation}  \label{9}
-(1+\theta p)e^{-(1+\theta p)}=-(1-\gamma)(1+\theta)e^{-\theta-1}.
\end{equation}
We see that $-(1+\theta p)$ is the Lambert-W function of real argument $-(1-\gamma)(1+\theta)e^{-\theta-1}$ .\\
Thus, we have
\begin{equation}  \label{10}
W\left(-(1-\gamma)(1+\theta)e^{-\theta-1}\right)=-(1+\theta p)
\end{equation}
Moreover, for any $\theta$, $\lambda$ and $\alpha >0 $ it is immediate that  $(1+p\theta)> 1 $, and it can also be checked that $(1-\gamma)(1+\theta)e^{-\theta-1} \in (-1/e,0)$ since $\gamma \in (0,1)$. Therefore, by taking into account the properties of the negative branch of the Lambert W function, we deduce the following. 
\begin{equation*}
W_{-1}\left(-(1-\gamma)(1+\theta)e^{-\theta-1}\right)=-(1+\theta p)
\end{equation*}
Again,solving for $x$ by  using $p=(1+\lambda x)^{\alpha}$, we get
\begin{equation}  \label{11}
x_{\gamma}=G^{-1}(\gamma)=\frac{1}{\lambda}\left(-\frac{W_{-1}\left(e^{-\theta -1} (\theta +1) (\gamma-1)\right)+1}{\theta }\right)^{1/\alpha}-\frac{1}{\lambda }
\end{equation} 
\begin{flushright}
$\Box$
\end{flushright}
Further the Median can be obtained by substituting $\gamma =\frac{1}{2}$ in (\ref{11}),Thus

\begin{equation}   \label{12}
\text{Median} (M_d)=G^{-1}(1/2)= \frac{1}{\lambda}\left(-\frac{W_{-1}\left(-\frac{e^{-\theta -1} (\theta +1) }{2}\right)+1}{\theta }\right)^{1/\alpha}-\frac{1}{\lambda } \\
\end{equation}

\subsection{Moments} 
An infinite sum representation is being used to represent $r^{th}$ moment, $\mu'_r = E[X^r]$,
and consequently the first four moments and variance for the EGL Distribution.

The $k^{th}$ raw moment of EGL random variable is given as \\
\begin{align} \label{13}
\nonumber \mathbb{E}(X^k)=&\frac{1}{1+\theta}\int_0^{\infty }\left(x^ke^{\theta \left(1-(1+x \lambda )^{\alpha }\right)} \alpha  \theta ^2 \lambda  (1+x \lambda )^{-1+2 \alpha }\right) dx \\ \nonumber
=& \frac{e^{\theta} \alpha  \theta ^2 \lambda }{1+\theta }\int _0^{\infty }x^ke^{-\theta  (1+x \lambda )^{\alpha }} (1+x \lambda )^{2 \alpha -1}dx \\
\mathbb{E}(X^k)=&\frac{e^{\theta} \alpha \theta ^2 \lambda }{1+\theta} I(k,\alpha ,\theta )
\end{align}
where $I(k,\alpha ,\theta)=\frac{1}{\alpha  \lambda ^{k+1}}\sum_{i=0}^k \binom{k}{i}(-1)^{k-i}\theta ^{-\frac{i}{\alpha }-2}\Gamma \left(\frac{i}{\alpha }+2,\theta \right)$ see appendix for detailed proof. \\
\\
Hence
\begin{equation}
\begin{aligned}
\mathbb{E}(X)= \frac{e^\theta}{\lambda (1+\theta)} &\left(-\Gamma (2,\theta)+\theta ^{-\frac{1}{\alpha }} \Gamma \left(2+\frac{1}{\alpha },\theta \right)\right)\\
\mathbb{E}(X^2)= \frac{e^{\theta }}{\lambda ^2(1+\theta)}& \left(\Gamma (2,\theta )-2 \theta ^{-\frac{1}{\alpha }} \Gamma \left(2+\frac{1}{\alpha },\theta \right)+\theta ^{-\frac{2}{\alpha }} \Gamma \left(2+\frac{2}{\alpha },\theta \right)\right) \\
\mathbb{E}(X^3)=\frac{e^{\theta }}{\lambda ^3(1+\theta )}& \left(-\Gamma (2,\theta )+3 \theta ^{-\frac{1}{\alpha }} \Gamma \left(2+\frac{1}{\alpha },\theta \right)-3 \theta ^{-\frac{2}{\alpha }} \Gamma \left(2+\frac{2}{\alpha },\theta \right)+\theta ^{-\frac{3}{\alpha }} \Gamma \left(2+\frac{3}{\alpha },\theta \right)\right) \\
\mathbb{E}(X^4)= \frac{e^{\theta }}{\lambda ^4(1+\theta)}& \left(\Gamma (2,\theta )-4 \theta ^{-\frac{1}{\alpha }} \Gamma \left(2+\frac{1}{\alpha },\theta \right)+6 \theta ^{-\frac{2}{\alpha }} \Gamma \left(2+\frac{2}{\alpha },\theta \right)-4 \theta ^{-\frac{3}{\alpha }} \Gamma \left(2+\frac{3}{\alpha },\theta \right)+\right.\\
&\left. \quad \theta ^{-\frac{4}{\alpha }} \Gamma \left(2+\frac{4}{\alpha },\theta \right)\right) 
\end{aligned}
\end{equation}

\noindent For lifetime models, it is also of interest to know moment of Future lifetime random variable $Y=X|X>t$ and its moments. Thus by using the Lemma the $k^{th}$ raw moment of random variable $Y$ is given as
\begin{equation}
\mathbb{E}(Y^k)=\mathbb{E}(X^k|X>t)= \frac{1}{\bar{G}(t)}\int _t^{\infty }u^k g(u)du
\end{equation}
Thus, substituting the value of $g(x)$ and $G(x)$ from equation () and () we get.
\begin{equation*}
\begin{aligned}
\mathbb{E}(X^k|X>t)=&\frac{\lambda \theta^2\alpha  e^{\theta (1+\lambda t)^{\alpha }}}{1+\theta (1+\lambda t)^{\alpha }}\int_t^{\infty}u^k(1+\lambda u)^{2\alpha -1}e^{-\theta (1+\lambda u)^{\alpha }}du \\
=& \frac{\lambda \theta ^2\alpha  e^{\theta (1+\lambda t)^{\alpha }}}{1+\theta (1+\lambda t)^{\alpha }} L(k,t,\alpha,\theta(1+\lambda t)^\alpha)
\end{aligned}
\end{equation*}
\noindent where $L(k,t,\alpha,\theta(1+\lambda t)^\alpha)=\frac{1}{\alpha \lambda ^{k +1}}\sum_{i=0}^k \binom{k}{i}(-1)^{k-i}\theta^{-\frac{i}{\alpha }-2}\Gamma \left(\frac{i}{\alpha }+2,\theta(1+\lambda t)^\alpha \right)$ see appendix for detailed proof. \\
The mean residual lifetime function is $\mathbb{E}(X|X>t)-t$.

\subsection{Entropy}
Let us now consider the R\'enyi entropy which represents a measure of uncertainty of a random variable and is defined as 
\begin{equation} \label{18} 
\mathfrak{J}(\zeta)=\frac{1}{1-\zeta} \log \left( \, \int\limits_{\mathbb{R}^+}f^\zeta(x) dx \right),  \quad \text{for} \quad \zeta >1 \quad \text{and} \quad \zeta \ne 1
\end{equation}
In our case 
\begin{equation*}
\int\limits_{\mathbb{R}^+} f^\zeta dx= \left(\frac{\alpha \lambda \theta ^2e^{\theta }}{1+\theta }\right)^{\zeta }\int _{\mathbb{R}^+}e^{-\zeta \theta (1+x \lambda )^{\alpha }}(1+x \lambda )^{\zeta (2\alpha -1)}dx
\end{equation*}
substituting $t=(1+\lambda x)^\alpha$ the above expression reduces to 
\begin{align} \label{19}
\nonumber \int\limits_{\mathbb{R}^+} f^\zeta dx=&\frac{(\alpha \lambda )^{\zeta -1}\theta ^{2\zeta }e^{\theta \zeta }}{(1+\theta )^{\zeta }}\int _1^{\infty }e^{-\gamma \theta t} t^{\frac{\gamma }{\alpha }(2\alpha -1)+\frac{1}{\alpha}-1}dt \\
=& \frac{(\alpha \lambda )^{\zeta -1}\theta ^{2\zeta }e^{\theta \zeta }}{(1+\theta )^{\zeta }} E_{\frac{-2 \zeta \alpha +\alpha +\zeta-1}{\alpha}}(\zeta \theta )
\end{align}
where $E_n(z)=\int_1^{\infty }e^{-zt} t^{-n} dt$ is known as exponential integral function. For more details \\
see http://functions.wolfram.com/06.34.02.0001.01. \\
Thus according to (\ref{18}) the R\'enyi entropy of EGL$(\theta,\lambda,\alpha)$ distribution given as 
\begin{equation} \label{20}
\mathfrak{J}(\zeta)=\frac{1}{1-\zeta}\log\left( \frac{(\alpha \lambda )^{\zeta -1}\theta ^{2\zeta }e^{\theta \zeta }}{(1+\theta )^{\zeta }} E_{\frac{-2 \zeta \alpha +\alpha +\zeta-1}{\alpha}}(\zeta \theta ) \right).
\end{equation}
Moreover, The Shannon entropy defined by $E[−\log(f(x))]$ is a special case derived from $\lim\limits_{\zeta \rightarrow 1}\mathfrak{J}(\zeta)$

\subsection{Order Statistics}
Let $X_1,X_2,\cdots,X_n$ be a random sample from the EGL$(\lambda, \theta, \alpha)$ distribution, and let $X_{i:n}$ denote the $i^{th}$ order statistic. Assuming $p=(1+\lambda x)^\alpha$, the pdf of the $i^{th}$ order statistic $X_{i:n}$ is given by (see [7])
\begin{equation*}
\begin{aligned}
g_{i:n}(x)=&\frac{n!}{(i-1)!(n-i)!}g(x)G^{i-1}(x)(1-G(x))^{n-i} \\
=&\frac{n!}{(i-1)!(n-i)!} \frac{\alpha \theta ^2\lambda }{(1+\theta )^{1+n-i}}e^{\theta (1-p)(1+n-i)}p^{2-\frac{1}{\alpha }}(1+\theta p)^{n-i}\left(1-\frac{e^{\theta (1-p)}(1+\theta p)}{(1+\theta )}\right)^{i-1} \\
=&\frac{n!}{(i-1)!(n-i)!}\frac{\alpha \theta ^2\lambda }{(1+\theta )^{1+n-i}}\sum _{j=0}^{i-1} \binom{i-1}{j} \frac{(-1)^j}{(1+\theta )^j}p^{2-\frac{1}{\alpha }}(1+\theta p)^{n-i+j}e^{\theta (1-p)(1+n-i+j)}
\end{aligned}
\end{equation*}
Substituting back the value of $p=(1+\lambda x)^\alpha$, we get
\begin{equation}
\begin{aligned}
\nonumber g_{i:n}(x)=\frac{n!}{(i-1)!(n-i)!}\frac{\alpha \theta ^2\lambda }{(1+\theta )^{1+n-i}}&\sum _{j=0}^{i-1} \binom{i-1}{j}\frac{(-1)^j}{(1+\theta )^j}(1+\lambda x)^{2\alpha -1}\left(1+\theta (1+\lambda x)^{\alpha }\right)^{n-i+j}\\
& \times \left(e^{\theta \left(1-(1+\lambda x)^{\alpha }\right)(1+n-i+j)}\right)
\end{aligned}
\end{equation}
Thus the moments of $X_{i:n}$ can be expressed as 
\begin{equation*}
\begin{aligned}
E\left(X_{i:n}^q\right)=&\frac{n!}{(i-1)!(n-i)!} \frac{\alpha \theta ^2\lambda }{(1+\theta )^{1+n-i}}\sum _{j=0}^{i-1} \sum _{k=0}^{n-i+j} \binom{i-1}{j}\binom{n-i+j}{k}\frac{(-1)^j}{(1+\theta )^j}e^{\theta (1+n-i+j)}\\
&\times \int _0^{\infty }x^q(1+\lambda x)^{\alpha (k+2)-1}e^{-\theta (1+n-i+j)(1+\lambda x)^{\alpha}}dx\\
=&\frac{n!}{(i-1)!(n-i)!} \frac{\alpha \theta ^2\lambda }{(1+\theta )^{1+n-i}}\sum _{j=0}^{i-1} \sum _{k=0}^{n-i+j} \binom{i-1}{j}\binom{n-i+j}{k}\frac{(-1)^j}{(1+\theta )^j}e^{\theta (1+n-i+j)}\\
&\times \frac{1}{\alpha \lambda ^{q+1}}\sum _{l=0}^q \binom{q}{l}(-1)^{q-l}(\theta (n+1-i+j))^{-\left(k+2+\frac{l}{\alpha }\right)}\Gamma \left(k+2+\frac{l}{\alpha},\theta (n+1-i+j)\right)
\end{aligned}
\end{equation*} 

\subsection{Limiting Distribution of Sample Minimum and Maximum}
Usually, we may be interested in the asymptotic behaviour of sample minima $X_{1:n}$ and/or sample maxima $X_{n:n}$.
Using, theorem 8.3.6 of Arnold et al., it follows that the asymptotic distribution of $X_{1:n}$ follows exponential whereas the $X_{1:n}$ follows  \\
Considering cdf given in equation() and strictly positive function $g(t)$ as $g(t)=\frac{1}{\lambda}(1+\lambda t)^{1-\alpha}, t>0 $.  It can be seen that \\

\begin{equation*}
\lim \limits_{t \to 0^+}\frac{F(tx)}{F(t)}=x ,\forall x >0
\end{equation*}
and 
\begin{align*}
\lim\limits_{t\rightarrow\infty}  \frac{1-F(t+xg(t))}{1-F(t)}=&
\lim\limits_{t\rightarrow\infty}  \frac{e^{-\theta (1+\lambda t+\lambda xg(t))^{\alpha}}\left(1+\theta(1+\lambda t+\lambda x g(t))^{\alpha}\right)}{1+\theta (1+\lambda t)^{\alpha}}\\
=&\lim\limits_{t\rightarrow\infty}  \frac{e^{-\theta (1+\lambda t)^{\alpha} \left(\left(\frac{\lambda x g(t)}{1+\lambda t}+1\right)^{\alpha}-1\right)}\left(\theta (1+\lambda+\lambda x g(t))^{\alpha}+1\right)}{1+\theta  (1+\lambda t)^{\alpha}}\\
=& \lim\limits_{t\rightarrow\infty} \frac{e^{-\theta (1+\lambda t)^{\alpha}\left[\left(\frac{\alpha \lambda x g(t)}{1+\lambda t}\right)+\cdots\right]}  \left(\theta  (1+\lambda t+\lambda x g(t))^{\alpha}+1\right)}{1+\theta  (1+\lambda t)^{\alpha}}\\
\end{align*}
\noindent substituting the value of $g(t)$ and taking limit, we obtain
\begin{align*}
\lim\limits_{t\rightarrow\infty}  \frac{1-F(t+xg(t))}{1-F(t)}=&\lim\limits_{t\rightarrow\infty}\frac{ \left(\theta  \left(1+\lambda t+x(\lambda t+1)^{1-\alpha}\right)^{\alpha}+1\right)}{1+\theta (\lambda t+1)^{\alpha}} e^{-\frac{{\alpha x} \theta  (\lambda  t+1)^{\alpha }+\cdots}{(\lambda t+1)^{\alpha }}}\\
=& e^{-\theta \alpha x}
\end{align*}
\indent so it follows from the theorem $1.6.2$ in Leaderbetter that there must be norming constants $a_n>0$ and $b_n$ such that $a_n=[g(\gamma_n)]^{-1}$ , $b_n=\gamma_n$ where $\gamma_n=F^{-1}(1-\frac{1}{n})$

\section{Maximum Likelihood Estimators}
In this section, we shall discuss the point and interval estimation of the parameters of EGL $ (\lambda,\theta,\alpha) $.
The log-likelihood function $\nonumber l(\Theta)$ of single observation (say $x_i$) for the vector of parameter $\Theta=(\theta, \lambda, \alpha)^\top$ is
\begin{equation*} 
l_n=n \log(\alpha)+n \log\lambda +2 n \log(\theta)-n \log(1+\theta)+\theta \sum _{i=1}^n \left(1-\left(1+x_i \lambda \right)^{\alpha}\right)+(2\alpha-1)\sum_{i=1}^n \log\left(1+x_i \lambda \right)
\end{equation*}
The associated score function is given by $U_n= \left(\frac{\partial l_n}{\partial \theta}, \frac{\partial l_n}{\partial \lambda}, \frac{\partial l_n}{\partial \alpha} \right)^\top$, where 

\begin{equation} \label{21}
\frac{\partial l_n}{\partial \theta}=\frac{2 n}{\theta }-\frac{n}{1+\theta}+\sum_{i=1}^n \left(1-(1+\lambda  x_i)^{\alpha}\right) \\
\end{equation}
\begin{equation} \label{22}
\frac{\partial l_n}{\partial \lambda}=\frac{n}{\lambda }+(2 \alpha-1)\sum_{i=1}^n \frac{x_i}{1+\lambda  x_i}-\theta\sum_{i=1}^n \alpha x_i (1+\lambda  x_i)^{\alpha-1}.
\end{equation}
\begin{equation} \label{}
\frac{\partial l_n}{\partial \alpha}=\frac{n}{\alpha}+2 \sum _{i=1}^n \log(1+\lambda  x_i)-\theta \sum _{i=1}^n (1+\lambda x_i)^{\alpha} \log(1+\lambda x_i) 
\end{equation}

As we know the expected value of score function equals  to zero, i.e.$E(U(\Theta))$, which implies\\
$\mathtt{E}\left(\sum_{i=1}^n \left(1-(1+\lambda  x_i)^{\alpha}\right)\right)= \frac{1}{\theta +1}-\frac{2}{\theta }$.

The total log-likelihood of the random sample $x=\left(x_1,\cdots, x_n\right)^\top$ of size $ n $ from $ X $ is given by $l_n= \sum\limits_{1}^{n}l^{(i)}$ and th total score function is given by $U_n=\sum\limits_{i=1}^{n}U^{(i)}$, where $\nonumber l^{(i)}$ is the log-likelihood of $i^{th}$ observation and $ U^{(i)} $ as given above. The maximum likelihood estimator $\hat{\Theta}$ of $\Theta$ is obtained by solving equations (\ref{21}) and (\ref{22}) numerically or this can also obtained easily by using nlm() function in R. The initial guess for the estimators were obtained from the inner region of 3D contour plot of log-likelihood function for a given sample. For example, in Figure (3), the contour plot of log-likelihood function for different $\theta$ and $\lambda$, the initial estimates were taken from interior. The associated Fisher information matrix is given by
%\begin{figure} \label{f3}
%\centering
%\includegraphics[width=0.9\textwidth]{contour}
%\caption{Contour plot of Log-likelihood surface for different values of $\theta$ and $\lambda$.}
%\end{figure}

\begin{equation} \label{23} 
K=K_n(\Theta)= n \left[ \begin{matrix}
\kappa_{\theta,\theta}&  \kappa_{\theta,\lambda} & \kappa_{\theta,\alpha} \\ 
\kappa_{\lambda,\theta}  & \kappa_{\lambda,\lambda} & \kappa_{\lambda,\alpha}\\
\kappa_{\alpha,\theta}  & \kappa_{\alpha,\lambda} & \kappa_{\alpha,\alpha}
\end{matrix} \right]
\end{equation}
where
\begin{equation} \label{24} 
\begin{aligned} 
\kappa_{\theta,\theta}&= \frac{n}{(\theta +1)^2}-\frac{2 n}{\theta ^2}\\
\kappa_{\lambda,\lambda}&= -\frac{n}{\lambda ^2}+\alpha \theta (1+\alpha)  \mathbb{E} \left( X^2 (\lambda X +1)^{\alpha -2} \right)+ (1-2 \alpha)  \mathbb{E} \left(\frac{X^2}{(\lambda  X+1)^2} \right)\\
\kappa_{\alpha,\alpha}&=-\frac{n}{\alpha ^2}-\theta \mathbb{E} \left( (\lambda  X+1)^{\alpha } \log (1+\lambda X)^2 \right)\\
\kappa_{\lambda,\theta}=\kappa_{\theta,\lambda}&= -\alpha \mathbb{E} \left( X (\lambda X+1)^{\alpha -1} \right) \\
\kappa_{\alpha,\theta}=\kappa_{\theta,\alpha}&= \mathbb{E} \left( (\lambda X + 1)^{\alpha} (-\log (\lambda X + 1)) \right)\\
\kappa_{\alpha,\lambda}=\kappa_{\lambda,\alpha}&= \mathbb{E} \left(-\theta X (\lambda X +1)^{\alpha -1}- \alpha \theta X (\lambda X+1)^{\alpha -1} \log (\lambda X+1)\right)+2  \mathbb{E} \left( \frac{X}{\lambda X+1} \right)\\
\end{aligned}
\end{equation}
The above expressions depend on some expectations which can be easily computed using numerical integration. Under the usual regularity conditions, the asymptotic distribution of
\begin{equation}
\sqrt{n}\left( \hat{\Theta}-\Theta\right) \, \,  \text{is} \, \,  N_2(0,K(\Theta)^{-1})\\ 
\end{equation}
where $\lim\limits_{n \rightarrow \infty}=K_n(\Theta)^{-1}=K(\Theta)^{-1}$. The asymptotic multivariate normal $N_2(0,K(\Theta)^{-1})$  distribution of $\hat{\Theta}$ can be used to construct approximate confidence  intervals. An asymptotic confidence interval with significance level $ \gamma $ for each
parameter $ \theta $,$ \alpha$ and $ \lambda $ is
\begin{equation}
\begin{aligned}
\text{ACI}\left(\theta,100(1-\gamma)\% \right)\equiv&(\hat{\theta}-z_{\gamma/2}\sqrt{\kappa_{(\theta,\theta)}},\hat{\theta}+z_{\gamma/2}\sqrt{\kappa_{(\theta,\theta)}}) \\
\text{ACI}\left(\lambda,100(1-\gamma)\%
\right)\equiv&(\hat{\lambda}-z_{\gamma/2}\sqrt{\kappa_{(\lambda,\lambda)}},\hat{\lambda}+z_{\gamma/2}\sqrt{\kappa_{(\lambda,\lambda)}})\\
\text{ACI}\left(\alpha,100(1-\gamma)\%
\right)\equiv&(\hat{\alpha}-z_{\gamma/2}\sqrt{\kappa_{(\alpha,\alpha)}},\hat{\alpha}+z_{\gamma/2}\sqrt{\kappa_{(\alpha,\alpha)}})\\
\end{aligned}
\end{equation}
where $z_{1-\gamma/2}$ denotes $1-{\gamma/2}$ quantile of standard normal random variable.
\section{Application to Real Datasets}
In this section, we illustrate, the applicability of Extended Generalized Lindley Distribution by considering two different datasets used by different researchers. We also fit Extended Generalized Lindley distribution, Lindley-Exponential Distribution proposed by Bhati et al.(2015), Power-Lindley distribution (2015), New Generalized Lindley Distribution(2013), Lindley Distribution(1958) and Exponential distribution. Namely\\

\noindent (i) Lindley-Exponential Distribution (L-E$(\theta,\lambda)$)
\begin{equation}
f(x)=\frac{ \theta ^2 \lambda e^{-\lambda  x} \left(1-e^{-\lambda x}\right)^{ \theta-1} \left(1-\log \left(1-e^{-\lambda x}\right)\right)}{(1+\theta)}  \quad \quad x,\theta,\lambda >0.
\end{equation}
(ii) Power-Lindley distribution (PL$(\alpha,\beta)$):
\begin{equation*}
f_1(x)=\frac{\alpha \beta^2}{1+\beta}(1+x^\alpha) x^{\alpha-1} e^{-\beta x^\alpha},  \quad \quad x,\alpha,\beta >0.
\end{equation*}

\noindent (iii) New Generalized Lindley distribution (NGLD($\alpha,\beta,\theta$)):
\begin{equation*}
f_2(x)=\frac{e^{-\theta x}}{1+\theta}\left(\frac{\theta^{\alpha+1}x^{\alpha-1}}{\Gamma(\alpha)}+\frac{\theta^{\beta}x^{\beta-1}}{\Gamma(\beta)}\right), \quad \quad x,\alpha,\theta,\beta >0.
\end{equation*}

\noindent(iv) Lindley Distribution (L$(\theta)$)
\begin{equation*}
f_3(x)=\frac{\theta^2}{1+\theta}(1+x) e^{-\theta x},  \quad \quad x,\alpha,\beta >0.
\end{equation*}

For each of these distributions, the parameters are estimated by using the maximum likelihood method, and for comparison we use Negative LogLikelihood values ($-LL$), the Akaike information criterion (AIC) and Bayesian information criterion (BIC) which are defined by $-2LL+2q$ and $-2LL+q\log(n)$, respectively, where $q$ is the number of parameters estimated and $n$ is the sample size. Further K-S statistic (Kolmogorov-Smirnov test)=$ \sup_x|F_n(x)-F(x)|$, where $F_n(x)=\frac{1}{n} \sum\limits_{i=1}^{n}\mathbb{I}_{x_i\le x}$ is empirical distribution function and $F(x)$ is cumulative distribution function is calculated and shown for all the datasets.

\subsection{Illustration 1: Application to bladder cancer patients}
We consider an uncensored data set corresponding to remission times (in months) of a random sample of 128 bladder cancer patients (Lee and Wang(2003)) as presented in Appendix A.1 in Table(9). The results for these data are presented in Table (\ref{t5}). We observe that the EGL distribution is a competitive distribution compared with other distributions. In fact, based on the values of the AIC and BIC criteria and K-S test statistic, we observe that the EGL distribution provides the best for these data among all the models considered. The probability density function and empirical distribution function are presented in figure (8) for all considered distributions for these data.

\begin{table}[htbp]
  \centering
  \caption{The estimates of parameters and goodness-of-fit statistics for Illustration 1.}
    \begin{tabular}{|l|l|c|c|c|c|} \hline
    Model & Parameter & -LL  & AIC   & BIC   & K-S statistic \\ \hline
    EGL & $ \hat{\lambda} $= 0.9360, $ \hat{\theta} $=0.5878,$ \hat{\alpha} $=0.6457 & \multicolumn{1}{c}{401.254} & \multicolumn{1}{|c|}{796.508} & \multicolumn{1}{|c|}{796.186} & \multicolumn{1}{|r|}{0.0387} \\
    L-E    & $ \hat{\lambda} $= 0.0962, $ \hat{\theta} $=1.229 & \multicolumn{1}{c}{401.78} & \multicolumn{1}{|c|}{807.564} & \multicolumn{1}{|c|}{807.780} & \multicolumn{1}{|r|}{0.0454} \\
    PL    & $\hat{\theta} $=0.385, $\hat{\beta}$=0.744 & \multicolumn{1}{|c|}{402.24} & \multicolumn{1}{|c|}{808.474} & \multicolumn{1}{|c|}{808.688} & \multicolumn{1}{|r|}{0.0446} \\
    L     & $\hat{\theta}$=0.196 & \multicolumn{1}{|c|}{419.52} & \multicolumn{1}{|c|}{841.040} & \multicolumn{1}{|c|}{843.892} & \multicolumn{1}{|r|}{0.0740} \\
    \multicolumn{1}{|l|}{NGLD} & $ \hat{\theta} $=0.180, $\hat{\alpha} $=4.679, $ \hat{\beta} $=1.324 & \multicolumn{1}{|c|}{412.75} & \multicolumn{1}{|c|}{831.501} & \multicolumn{1}{|c|}{840.057} & \multicolumn{1}{|r|}{0.1160} \\ \hline
    \end{tabular}%
  \label{t5}%
\end{table}%

\begin{figure} \label{f8}
\centering
%\textit{}\includegraphics[width=1\textwidth]{cdf_hist_128}
\includegraphics[width=1\textwidth]{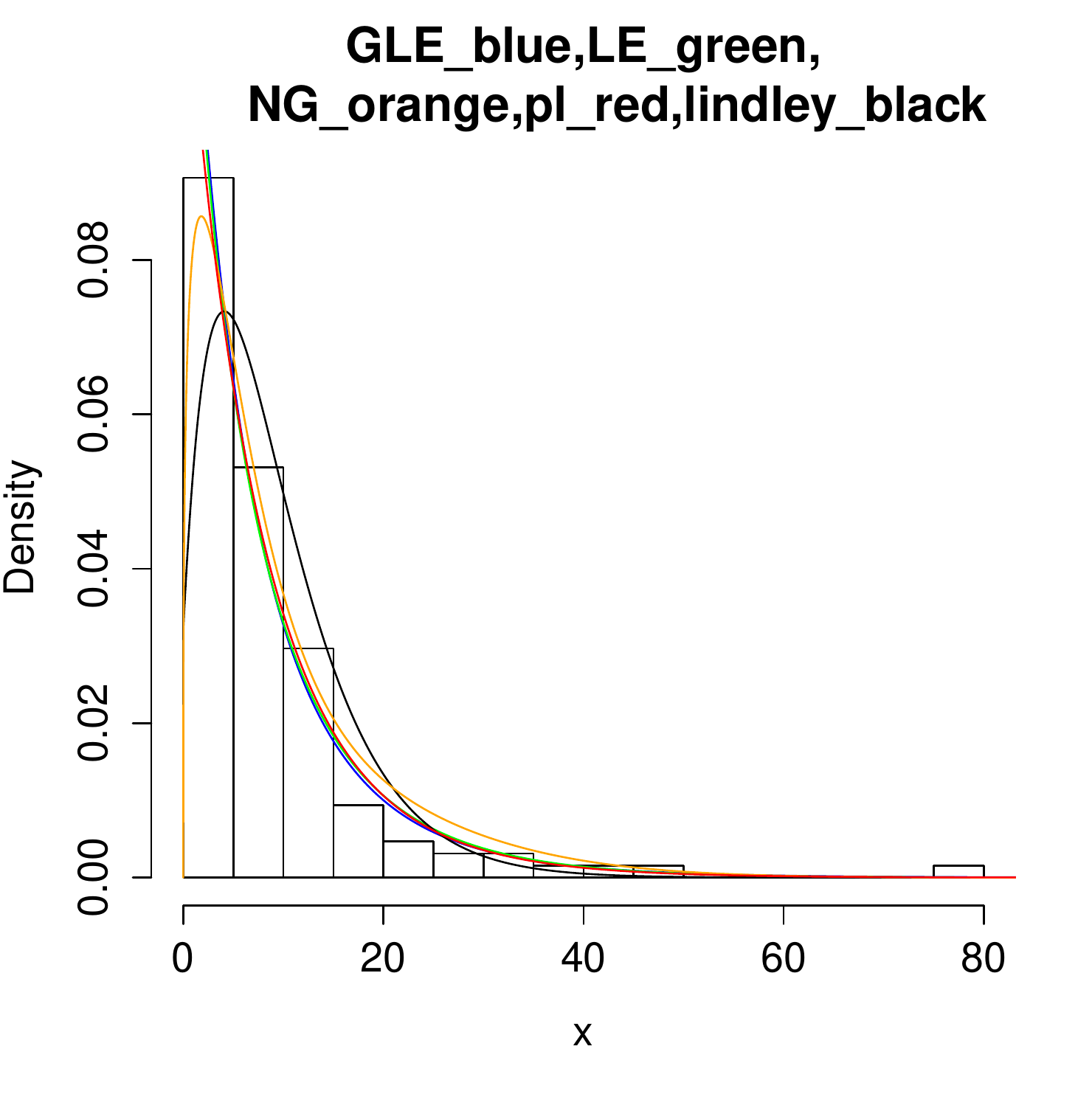}
\caption{PDF plot for various values of $\lambda$ and $\theta$.}
\end{figure}

\subsection{Illustration 2: Application to waiting times in a queue}
As second example, we consider 100 observations on waiting time (in minutes) before the customer received service in a bank (see Ghitany et al.(2008)). The data sets are presented in appendix A.2 as Table (10). The results for these data are presented in Table(\ref{t6}). From these results we can observe that EGL distribution provide smallest K-S test statistics values as compare to Lindley-Exponential, Power lindley, new generalized Lindley distribution, Lindley and Exponential and hence best fits the data among all the models considered. The results are presented in Table (\ref{t6}) and probability density function and empirical distribution function are shown in figure (9).
\begin{table}[htbp]
  \centering
  \caption{The estimates of parameters and goodness-of-fit statistics for Illustration 2.}
    \begin{tabular}{|l|l|r|r|r|r|} \hline
    Model & Parameter & LL    & AIC   & BIC   & K-S \\ \hline
    EGL & $ \hat{\lambda} $= 1.803, $ \hat{\theta} $=0.093,$ \hat{\alpha} $=1.046 & \multicolumn{1}{c}{318.066} & \multicolumn{1}{|c|}{642.132} & \multicolumn{1}{|c|}{642.132} & \multicolumn{1}{|r|}{0.0384} \\
    L-E    & $\hat{\theta}$=2.650, $\hat{\lambda}$=0.1520 & 317.005 & 638.01 & 638.1337  & 0.0360 \\
    PL    & $\hat{\theta}$=0.1530;$\hat{\beta}$=1.0832 & 318.319 & 640.64  & 640.64 & 0.0520 \\
    L & $\hat{\theta}$=0.187 & 319.00   & 640.00   & 640.00   & 0.0680 \\
    E & $\hat{\theta}$=0.101 & 329.00   & 660.00   & 660.00   & 0.1624 \\
    NGLD  & $\hat{\theta}$= 0.2033; $\hat{\beta}$=2.008; $\hat{\alpha}$=2.008 & 317.3 & 640.60  & 640.60 & 0.0425 \\ \hline
    \end{tabular}%
  \label{t6}%
\end{table}%
\begin{figure} \label{f9}
\centering
\includegraphics[width=1\textwidth]{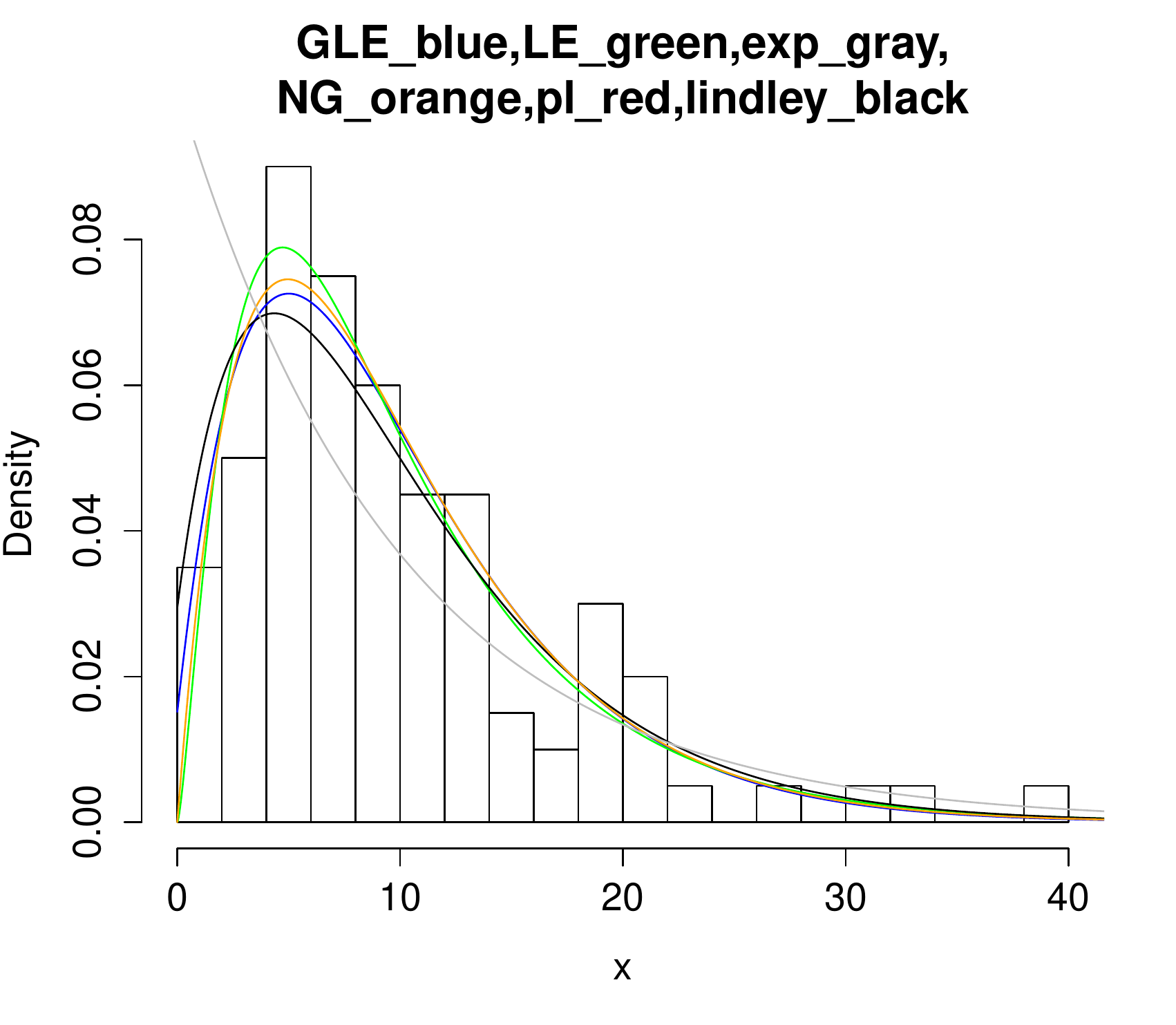}
\caption{PDF plot for various values of $\lambda$ and $\theta$.}
\end{figure}

\section*{Conclusion}
We have proposed a new three parameters class of distributions called Extended Generalized Lindley (EGL) distribution generated by Lindley distribution  which possess increasing, decreasing or upside down hazard function for different choices of the parameters. We have derived important properties of the EGL distribution like moments, entropy, asymptotic distribution of sample maximum and sample Minimum. Maximum likelihood of the parameters are obtained which can be used to get asymptotic confidence intervals. We have also illustrated the application of EGL distribution to two real data sets used by researchers earlier and compare it with other popular models. Further the stress-strength analysis were carried out and compared with that of Lindley distribution. Our application to real data set indicate that EGL distribution performs satisfactorily or better than its competitors and can be recommended for lifetime modelling the encountered in engineering, medical science, biological science and other applied sciences.

\section*{Appendix}

\noindent \textbf{A.1- Dataset used in Illustration 1:}

\begin{small}
\begin{table}[htbp]
  \centering
  \caption{The remission times (in months) of bladder cancer patients}
    \begin{tabular}{|r|r|r|r|r|r|r|r|r|r|r|r|r|} \hline
    0.08  & 2.09  & 3.48  & 4.87  & 6.94  & 8.66  & 13.11 & 23.63 & 0.2   & 2.23  & 0.26  & 0.31  & 0.73 \\ \hline
    0.52  & 4.98  & 6.97  & 9.02  & 13.29 & 0.4   & 2.26  & 3.57  & 5.06  & 7.09  & 11.98 & 4.51  & 2.07 \\ \hline
    0.22  & 13.8  & 25.74 & 0.5   & 2.46  & 3.64  & 5.09  & 7.26  & 9.47  & 14.24 & 19.13 & 6.54  & 3.36 \\ \hline
    0.82  & 0.51  & 2.54  & 3.7   & 5.17  & 7.28  & 9.74  & 14.76 & 26.31 & 0.81  & 1.76  & 8.53  & 6.93 \\ \hline
    0.62  & 3.82  & 5.32  & 7.32  & 10.06 & 14.77 & 32.15 & 2.64  & 3.88  & 5.32  & 3.25  & 12.03 & 8.65 \\ \hline
    0.39  & 10.34 & 14.83 & 34.26 & 0.9   & 2.69  & 4.18  & 5.34  & 7.59  & 10.66 & 4.5   & 20.28 & 12.63 \\ \hline
    0.96  & 36.66 & 1.05  & 2.69  & 4.23  & 5.41  & 7.62  & 10.75 & 16.62 & 43.01 & 6.25  & 2.02  & 22.69 \\ \hline
    0.19  & 2.75  & 4.26  & 5.41  & 7.63  & 17.12 & 46.12 & 1.26  & 2.83  & 4.33  & 8.37  & 3.36  & 5.49 \\ \hline
    0.66  & 11.25 & 17.14 & 79.05 & 1.35  & 2.87  & 5.62  & 7.87  & 11.64 & 17.36 & 12.02 & 6.76  &  \\ \hline
    0.4   & 3.02  & 4.34  & 5.71  & 7.93  & 11.79 & 18.1  & 1.46  & 4.4   & 5.85  & 2.02  & 12.07 &  \\ \hline
    \end{tabular}%
\end{table}
\label{t9}
\end{small}

\noindent \textbf{A.2- Dataset used in Illustration 2:}
\begin{small}
\begin{table}[htbp]
  \centering
  \caption{Waiting times (min.) of 100 bank customers}
    \begin{tabular}{|r|r|r|r|r|r|r|r|r|r|} \hline
    0.8   & 0.8   & 1.3   & 1.5   & 1.8   & 1.9   & 1.9   & 2.1   & 2.6   & 2.7 \\
    2.9   & 3.1   & 3.2   & 3.3   & 3.5   & 3.6   & 4     & 4.1   & 4.2   & 4.2 \\
    4.3   & 4.3   & 4.4   & 4.4   & 4.6   & 4.7   & 4.7   & 4.8   & 4.9   & 4.9 \\
    5.0   & 5.3   & 5.5   & 5.7   & 5.7   & 6.1   & 6.2   & 6.2   & 6.2   & 6.3 \\
    6.7   & 6.9   & 7.1   & 7.1   & 7.1   & 7.1   & 7.4   & 7.6   & 7.7   & 8 \\
    8.2   & 8.6   & 8.6   & 8.6   & 8.8   & 8.8   & 8.9   & 8.9   & 9.5   & 9.6 \\
    9.7   & 9.8   & 10.7  & 10.9  & 11.0  & 11.0  & 11.1  & 11.2  & 11.2  & 11.5 \\
    11.9  & 12.4  & 12.5  & 12.9  & 13.0    & 13.1  & 13.3  & 13.6  & 13.7  & 13.9 \\
    14.1  & 15.4  & 15.4  & 17.3  & 17.3  & 18.1  & 18.2  & 18.4  & 18.9  & 19.0 \\
    19.9  & 20.6  & 21.3  & 21.4  & 21.9  & 23    & 27    & 31.6  & 33.1  & 38.5 \\ \hline
    \end{tabular}%
\end{table}
\label{t10}
\end{small}


\begin{thebibliography}{}
\bibitem{r1} Adamidis K., and Loukas S.(1998) A lifetime distribution with decreasing failure rate, \textit{Statistics and Probability Letters}, \textbf{(39)}, 35-42.

\bibitem{r2} Al-Mutairi, D.K., Ghitany, M.E. and Kundu, D.(2013) Inference on stress- strength reliability from Lindley distribution,\textit{Communications in Statistics - Theory and Methods}, \textbf{(42)}, 1443-1463.

\bibitem{r3} Arnold B.C., Balakrishnan N. and Nagaraja H.N.(2013) A First Course in Order Statistics, \textit{Wiley, New York}.

\bibitem{r4} Bakouch H. S., Al-Zahrani B. M., Al-Shomrani A. A., Marchi V. A., and Louzada F.(2012) An extended Lindley distribution,\textit{Journal of the Korean Statistical Society}, \textbf{41}, 75-85.

\bibitem{r5} Bhati D. Malik A. and Vaman H.J.(2015) On Lindley-Exponential Distribution: Application and Properties, Metron, 73(3), 335-357. 
\bibitem{r6} Elbatal I., Merovci F., and Elgarhy M.(2013) A new generalized Lindley distribution, \textit{Mathematical Theory and Modeling}, \textbf{3(13)}.

\bibitem{r7} Ghitany M. E., Alqallaf F., Al-Mutairi D. K., and Husain H. A.(2011) A two-parameter weighted Lindley distribution and its applications to survival data, \textit{Mathematics and Computers in Simulation}, \textbf{81(6)}, 1190-1201.

\bibitem{r8} Ghitany M. E., Atieh B., and Nadarajah, S.(2008) Lindley distribution and its application, \textit{Mathematics and Computers in Simulation}, \textbf{78}, 493-506.

\bibitem{r9} Ghitany M. E., Al-Mutairi D. K. and Aboukhamseen S. M.(2015) Estimation of the reliability of a stress-strength system from power lindley distributions, \textit{Communications in Statistics-Simulation and Computation}, \textbf{44(1)}, 118-136.

\bibitem{r10} Glaser R.E.(1980) Bathtub and related failure rate characterizations, \textit{Journal of American Statistical Association}, \textbf{75} ,667–672.

\bibitem{r11} G\'omez E. D., Sordo M. A., and Calder\'in E. O.(2014) The Log–Lindley distribution as an alternative to the beta regression model with applications in insurance, \textit{Insurance: Mathematics and Economics}, \textbf{54}, 49-57.

\bibitem{r12} Hassan M.K.(2014) On the Convolution of Lindley Distribution,\textit{Columbia International Publishing Contemporary Mathematics and Statistics}, \textbf{2(1)}, 47-54.

\bibitem{r13} Jo\'rda P.(2010) Computer generation of random variables with Lindley or Poisson-Lindley distribution via the Lambert W function, \textit{Mathematics and Computers in Simulation}, \textbf{81}, 851-859.

\bibitem{r14} Krishnamoorthy, K., Mukherjee, S., Guo, H.(2007) Inference on reliability in two parameter exponential stress-strength model, \textit{Metrika}, \textbf{65}, 261-273.

\bibitem{r15} Kundu, D., Gupta, R. D.(2005) Estimation of $R= P(Y < X)$ for generalized exponential distributions, \textit{Metrika}, \textbf{61}, 291–380.

\bibitem{r16} Kundu, D., Gupta, R. D.(2006) Estimation of $R= P(Y < X)$ for Weibull distributions,
\textit{IEEE Trans. Reliab.}, \textbf{55}, 270–280.

\bibitem{r17} Kundu, D., Raqab, M. Z.(2009) Estimation of $R=P(Y < X)$ for three-parameter Weibull distribution, \textit{Statistics Probability Letters}, \textbf{79},1839–1846.

\bibitem{r18} Leadbetter M. R., Lindgren G., Rootzén H.(1983) Extremes and Related Properties of Random Sequences and Processes, \textit{Springer Statist. Ser., Springer}, Berlin.

\bibitem{r19} Lee E.T., and Wang J.W.(2003), Statistical methods for survival data analysis, \textit{John Wiley \& Sons}, inc., Hoboken, New Jersey, 3ed.

\bibitem{r20} Lindley D. V.(1958) Fiducial distributions and Bayes’ theorem, \textit{Journal of the Royal Statistical Society}, \textit{Series B (Methodological)}, 102-107.

\bibitem{r21} Nair, N. Unnikrishnan, Sankaran, P.G., Balakrishnan, N.(2013) Quantile-Based Reliability Analysis,
\textit{Springer}.

\bibitem{r22} Raqab, M. Z., Kundu, D.(2005) Comparison of different estimators of $P(Y < X)$ for a scaled Burr type X distribution, \textit{Communication in Statistics: Simulation and Computation}, \textbf{34}, 465-483.

\bibitem{r23} Risti\'c M. M. and Balakrishnan N.(2012) The gamma exponentiated exponential distribution, \textit{Journal of Statistical Computation and Simulation}, \textbf{82(8)}, 1191-1206.

\bibitem{r24} Shanker R., Sharma S., and Shanker R.(2013) A Two-Parameter Lindley Distribution for Modeling Waiting and Survival Times Data, \textit{Applied Mathematics}, \textbf{4}, 363-368.

\bibitem{r25} Zakerzadeh H. and Mahmoudi E.(1998) A new two parameter lifetime distribution: model and properties, arXiv:1204.4248 v1 [stat.CO], 2012.
\end{thebibliography}
\end{document}